\documentclass[a4paper,10pt]{article}

\usepackage{amssymb,amsmath,dsfont}

\newcommand{\bs}{\mathbf{s}}

\newcommand{\clS}{\mathcal{S}}

\newcommand{\hcf}{{\rm h.c.f.}}

\newcommand{\adj}{{\rm adj}}

\renewcommand{\b}[1]{{\bf #1}}

\renewcommand{\a}{\mathbf{a}}

\newcommand{\Z}{\mathbb{Z}}

\newcommand{\C}{\mathbb{C}}

\newcommand{\N}{\mathbb{N}}

\newcommand{\R}{\mathbb{R}}

\newcommand{\ep}{\varepsilon}

\newcommand{\x}{\mathbf{x}}

\newcommand{\beql}[1]{\begin{equation}\label{#1}}

\newcommand{\eeq}{\end{equation}}

\newcommand{\andd}{\;\;\;\mbox{and}\;\;\;}
\newcommand{\modd}[1]{\; ( \text{mod} \; #1)}

\newtheorem{theorem}{Theorem}
\numberwithin{equation}{section}
\newtheorem{lemma}{Lemma}

\begin{document}
\title{Small Solutions of Quadratic Congruences, and Character Sums
  with Binary Quadratic Forms}
\author{D.R. Heath-Brown\\Mathematical Institute, Oxford}
\date{}
\maketitle
\section{Introduction}

Let $Q(\x)=Q(x_{1},\ldots,x_{n})\in\Z [x_{1},\ldots,x_{n}]$
be a quadratic form.  This paper, which may be seen as a continuation
of the author's earlier work \cite{ssqc1}, \cite{ssqc2} seeks to
understand the
smallest solution of the congruence $Q(\x)\equiv 0\modd{q}$ in
non-zero integers $\x$.
Thus we shall set
\[m(Q;q):=\min\{||\x||:\x\in\Z^n-\{\mathbf{0}\},\,Q(\x)\equiv
0\modd{q}\}\]
where $||\x||$ denotes the Euclidean norm, and ask (in the first instance) about
\[B_n(q):=\max_Q m(Q;q),\]
where the maximum is taken over all integral quadratic forms in $n$
variables. (This definition differs slightly from that used in
\cite{ssqc1} and \cite{ssqc2}.)  The interested reader may
refer to Baker \cite[Chapter 9]{baker} for an account of this problem and its 
applications.

It is trivial that $B_n(q)$ is non-increasing as a function of
$n$. When $q$ is square-free it is easy to see that $B_n(q)\ge q$
for $n=1$ or 2.  
Moreover the form $Q(\x)=x_1^2+\ldots+x_n^2$ has $m(Q;q)\ge
q^{1/2}$ so that $B_n(q)\ge q^{1/2}$ for every $q$ and every $n$. When
$n=3$ and $q$ is square-free one has
\[B_3(q)\ge m(Q;q)\ge q^{2/3}+O(q^{1/3})\]
for a suitable singular form
\beql{ex}
Q(x_1,x_2,x_3)=(x_1-bx_2)^2-a(x_2-bx_3)^2.
\eeq
(Details for the case in which $q$ is prime are given in 
\cite[Theorem 3]{ssqc1} but the argument readily extends to any
square-free $q$.) It is reasonable to
conjecture that such lower bounds represent the true order of
magnitude for $B_n(q)$ in general, so that one would have
\[B_n(q)\ll_{\ep}\left\{\begin{array}{cc}
q^{2/3+\ep}, & n=3\\ q^{1/2+\ep},& n\ge 4
\end{array}\right.\]
for any fixed $\ep>0$ (uniformly in $n$, by the non-increasing property).

A basic upper bound for $B_n(q)$ was provided by Schinzel, Schlickewei 
and Schmidt \cite{SSS}, who showed that
\[B_n(q)\ll\left\{\begin{array}{cc}
q^{1/2+1/(2n)}, & n \mbox{ odd,}\\ q^{1/2+1/(2n-2)},& n \mbox{ even}.
\end{array}\right.\]
In particular one sees that $q^{2/3}$ is the true order of magnitude of $B_3(q)$,
at least when $q$ is square-free.

For $n\ge 4$ and any $\ep>0$ one has
\[B_n(q)\ll_{\ep} q^{1/2+\ep}\]
if $q$ has at most 2 prime factors, (see the author \cite[Theorem
1]{ssqc2}); that
\[B_4(q)\ll_{\ep} q^{5/8+\ep}\]
(see \cite[Theorem 2]{ssqc2}); and that
\[B_n(q)\ll_{\ep,n} q^{1/2+3/n^2+\ep}\]
for every even $n\ge 2$ (see \cite[Theorem 3]{ssqc2}). Indeed a number
of other such bounds are possible.

It might appear from the above discussion that our question is
completely resolved for $n=3$, but wide open for $n\ge 4$.  None the
less, the
main goal of this paper is a further exploration of the situation for
$n=3$ (!) It will be observed that the example (\ref{ex}) is a
singular form.  It turns out that one can do better if one restricts
attention to ternary forms which are nonsingular modulo $q$. Before
stating our result we should make two simple observations.  Firstly,
if $q=q_0^2q_1$ and $q_1\mid Q(\x)$, then $q\mid Q(q_0\x)$.  It
follows that $m(Q;q)\le q_0m(Q;q_1)$.  In particular, if we have proved
that $B_n(q)\ll q^{\theta}$ for all square-free $q$, for some exponent
$\theta\ge \tfrac12$, then we may deduce that $B_n(q)\ll q^{\theta}$
for every $q$.  Secondly, if $q=2q_1$ and $q_1\mid Q(\x)$, then 
$q\mid Q(2\x)$.  It follows in this case that $m(Q;q)\le 2m(Q;q_1)$.
Once again, if we have proved
that $B_n(q)\ll q^{\theta}$ for all odd square-free $q$, for some exponent
$\theta$, then we may deduce that $B_n(q)\ll q^{\theta}$
for every square-free $q$. These observations allow us to focus on odd
square-free $q$. Indeed we shall assume without further comment
throughout the remainder of this paper that $q$ is odd and
square-free. In this situation we can represent $Q(\x)$ modulo $q$
via a symmetric integer matrix, which we also denote by $Q$, by abuse
of notation.  

We now define
\[B_3^*(q):=\max_Q m(Q;q),\]
where the maximum is taken over all integral ternary quadratic forms
$Q$ with $(\det(Q),q)=1$.
This notation allows us to state our principal result.
\begin{theorem}\label{main}
Let $q\in\N$ be odd and square-free, and let $\ep>0$ be given.
Then
\[B_3^*(q)\ll_{\ep} q^{5/8+\ep}.\]
\end{theorem}
So we see that we can go below the exponent $2/3$ which is
the limiting exponent for $B_3(q)$. We now have the same exponent
$5/8$ for (non-singular) forms in 3 variables as we previously had for
4 variables. (However it is explained in \cite{ssqc2} that one can
reduce the exponent to $13/21$ with more work, in the 4 variable case.)

It now seems that one should conjecture a bound
\beql{c1}
B_3^*(q)\ll_{\ep} q^{1/2+\ep}.
\eeq

The proof of Theorem \ref{main} proceeds by reducing the problem to a
second question, which we now explain. If $Q$ is a 
quadratic form in $n\ge 2$ variables we write
\[\widehat{m}(Q;q):=\min\{||\x||:\x\in\Z^n-\{\mathbf{0}\},\,\exists
t\in\Z,\,Q(\x)\equiv t^2\modd{q}\}\]
and 
\[\widehat{B}_n(q):=\max_Q \widehat{m}(Q;q),\]
where the maximum is taken over all integral quadratic forms in $n$
variables such that $(\det(Q),q)=1$.

We then have the following result.
\begin{lemma}\label{key}
Let $q\in\N$ be odd and square-free.  Then if $Q$ is a ternary quadratic form
with $(\det(Q),q)=1$ we have
\[m(Q;q)\ll q^{1/2}\widehat{m}(-Q^{\adj};q)^{1/2},\]
where $Q^{\adj}$ is the adjoint matrix for $Q$.  In particular one has
\[B_3^*(q)\ll q^{1/2}\widehat{B}_3(q)^{1/2}.\]
\end{lemma}
This naturally leads us to speculate about the size of $\widehat{B}_3(q)$,
and the natural conjecture is that
\beql{c1a}
\widehat{B}_3(q)\ll_{\ep} q^{\ep}
\eeq
for any fixed $\ep>0$.  Of course Lemma \ref{key} immediately shows
that this latter conjecture implies (\ref{c1}).

If $q$ is odd and square-free there is a real character 
\[\chi_d(m)=\left(\frac{m}{d}\right)\]
for each divisor $d$ of $m$,
and the congruence $Q(\x)\equiv t^2\modd{q}$ will have a solution $t$
if and only if
\[\sum_{d\mid q} \chi_d(Q(\x))>0.\]
We can therefore attempt to show that $\widehat{B}_3(q)$ is small by
investigating the character sums
\[S(\chi,B,Q):=\sum_{||\x||\le B}\chi(Q(\x))\]
for primitive characters $\chi$ to modulus $d>1$.  
If we can show that 
\beql{c2}
S(\chi,B,Q)\ll B^{3-\delta}
\eeq
for some fixed
$\delta>0$, for every primitive $\chi$ to modulus $d\ge 2$, 
then we will be able to
deduce that $\widehat{m}_3(Q;q)\le B$, since one has
$S(1,B,Q)\gg B^3$ for the trivial character. 

It seems plausible that (\ref{c2}) should hold for $B\ge q^{\eta}$,
for any fixed $\eta>0$, and with $\delta=\delta(\eta)>0$.  This would
suffice for (\ref{c1a}), and hence also for (\ref{c1}). 

One standard
procedure to estimate sums such as $S(\chi,B,Q)$ is to complete the
sum and use bounds of Weil--Deligne type.  It is very instructive to
carry this out in detail.  What one finds, if $d=q$ for example, 
is essentially that 
\[S(\chi,B,Q)\ll \frac{B^3}{q^3}\sum_{\b{y}\in\Z^3}W\left(\frac{q}{B}\b{y}\right)
S(\b{y}),\]
where $W\ll 1$ is a suitable weight function and
\[S(\b{y}):=\sum_{\x\modd{q}}e_q(\b{y}.\x)\chi(Q(\x)).\]
These complete sums can be computed explicitly.  Taking $q$ to be prime
for simplicity, and assuming that $q\nmid \det(Q)$, one finds that
$S(\b{y})$ is of order $q^2$ when $q\mid Q^{\adj}(\b{y})$, and of
order $q$ otherwise. This may be something of a surprise, since one
typically expects complete sums in $n$ variables to have size around
$q^{n/2}$. As a result this analysis leads to a bound which one may
think of as
\[S(\chi,B,Q)\ll_{\ep}q^{\ep}(q+B^3q^{-1}\#\{\b{y}\ll q/B:\, q\mid
Q^{\adj}(\b{y})\}).\]
Since we have estimated $m(Q;q)$ in terms of sums $S(\chi,B,Q^{\adj})$
it is apparent that the above analysis ultimately connects small
solutions of $q\mid Q(\x)$ with small solutions of $q\mid Q(\b{y})$.
In fact the argument is not completely circular, and one can show in
this way that $B_3^*(q)\ll_{\ep} q^{2/3+\ep}$, at least when $q$ is
prime. Alternatively one can provide an upper bound for 
\[\#\{\b{y}\ll q/B:\, q\mid Q^{\adj}(\b{y}),\,\b{y} \mbox{ primitive}\}\]
by using $O((q/B)^{3/2})$ plane slices of the type $\b{a}.\b{y}=0$.
Each such slice produces a binary quadratic form of rank 1 or 2, which
will have $O(1)$ primitive zeros modulo $q$, under the assumption that
$q$ is prime.  In this way one finds that
\[\#\{\b{y}\ll q/B:\, q\mid Q^{\adj}(\b{y})\}\ll (q/B)^{3/2},\]
and hence
\[S(\chi,B,Q)\ll_{\ep}q^{\ep}(q+B^{3/2}q^{1/2}).\]
We therefore have a non-trivial bound for $B\ge
q^{1/3+\delta}$. Unfortunately this merely yields
$\widehat{B}_3(q)\ll_{\ep}q^{1/3+\ep}$ and hence $B_3^*(q)\ll q^{2/3+\ep}$.

We have been unable to obtain a non-trivial bound for $S(\chi,B,Q)$
when $B\le q^{1/3}$. However, if one replaces $Q$ by a binary form one
can do better. Indeed the following result of Chang \cite[Theorem 11]{ch}
is the main inspiration for this paper.
\begin{theorem}\label{thch}
{\bf (Chang).}
For any $\ep>0$ there is a corresponding $\delta>0$ such that
\[\left|\sum_{X'<x\le X+X'}\sum_{Y'<y\le Y+Y'}\chi(x^2+axy+by^2)\right|<
p^{-\delta}XY\]
for any non-trivial character $\chi$ modulo $p$, any 
$X,Y>p^{1/4+\ep}$, and any integers $a,b$ with 
$a^2\not\equiv 4b\modd{p}$.
\end{theorem}
This improves on the corresponding results of Burgess \cite{bq1} and
\cite {bq2}, which were non-trivial only for $X,Y>p^{1/3+\ep}$.

The proof of Chang's result
crucially uses the fact that a binary quadratic form over
$\mathbb{F}_p$ factors over $\mathbb{F}_{p^2}$, and of course this
limits the approach to the case $n=2$.  Since we are interested in
composite $q$ we will require a variant of Theorem \ref{thch}. The
argument in \cite{ch} splits into two rather different cases, one in 
which the form factors over $\mathbb{F}_p$, and one in which it does
not. In order to handle composite $q$ we need to devise a treatment
which handles both cases in the same way.  Our result is the following.

\begin{theorem}\label{burgthm}
Let $\ep>0$ and an integer $r\ge 3$ be given, and suppose that 
$C\subset \R^2$ is a convex
set contained in a disc $\{\x\in\R^2:\,||\x-\x_0||\le R\}$.  
Let $q\ge 2$ be odd and
square-free, and let $\chi$ be a primitive character to modulus
$q$.  Then if $Q(x,y)$ is a binary quadratic form with $(\det(Q),q)=1$ we have
\[\sum_{(x,y)\in
  C}\chi(Q(x,y))\ll_{\ep,r}R^{2-1/r}q^{(r+2)/(4r^2)+\ep}\]
for $q^{1/4+1/2r}\le R\le q^{5/12+1/2r}.$
\end{theorem}

For comparison we observe that the standard Burgess bound
\cite[Theorem 2]{burg} yields
\[\sum_{x,y\le R}\chi(xy)\ll_{\ep,r}R^{2-2/r}q^{(r+1)/(2r^2)+\ep},\]
relative to which our theorem has a loss of $(Rq^{-1/4})^{1/r}$.
In Section \ref{3}
we will apply Theorem \ref{burgthm} with $C=\{\x\in\R^2:\,||\x||\le R\}$
and $R=q^{1/4+\delta}$. Taking $r>(2\delta)^{-1}$
we will be able to deduce that $\widehat{B}_2(q)\ll_{\ep} q^{1/4+\delta}$.  We then
go on to conclude that
 $\widehat{B}_3(q)\ll_{\ep}q^{1/4+\ep}$ and hence, via Lemma \ref{key}, that
$B_3^*(q)\ll q^{5/8+\ep}$. 

Before embarking on the proofs we need to mention one point of
notation.  We shall follow the common convention that the small
positive number $\ep$ will be allowed to change between appearances,
allowing us to write $q^{\ep}\log q\ll_{\ep}q^{\ep}$, for example.

\section{Proof of Lemma \ref{key}}

Suppose that $-Q^{\adj}(\a)\equiv t^2\modd{q}$ with  
$||\a||= \widehat{m}(-Q^{\adj};q)$ and $\a\not=\mathbf{0}$.  
Write $\a=\alpha\a_0$ with
$\alpha\in\N$ and $\a_0$ primitive, and let
\[\Lambda:=\{\x\in\Z^3:\a_0.\x=0\}\]
This will be a 2-dimensional lattice of determinant $||\a_0||$.
Let $\x_1$ be the shortest non-zero vector in $\Lambda$, and $\x_2$
the shortest vector non-proportional to $\x_1$.  Then $\x_1$ and
$\x_2$ form a basis for $\Lambda$, and we have 
\beql{12a}
||\x_1||.||\x_2||\ll ||\a_0||
\eeq
and
\[\x_1\wedge\x_2=\pm\a_0.\]
We proceed to write
$R(u,v):=Q(u\x_1+v\x_2)$, so that $R$ is a binary quadratic form. We
then have
\[\det(R)=\det(Q(u\x_1+v\x_2))=Q^{\adj}(\x_1\wedge\x_2)=Q^{\adj}(\a_0)\]
as an identity, so that 
\[-\alpha^2\det(R)=-Q^{\adj}(\a)\equiv t^2\modd{q}. \]
Let $(q,\alpha)=q_0$ and $q=q_0q_1$.  It follows that
$R$ factors over $\mathbb{F}_p$ for every prime factor $p$ of $q_1$.  We
may then use the Chinese Remainder Theorem to
write $R(u,v)\equiv L_1(u,v)L_2(u,v)\modd{q_1}$ for certain
integral linear forms $L_1$ and $L_2$.

Our strategy is now to find a short vector $\x\in\Lambda$ such that
$q_1\mid L_1(u,v)$. We will then automatically have $q_1\mid R(u,v)$ and
hence $q_1\mid Q(u\x_1+v\x_2)$. This will
produce $q\mid Q(\x)$ with $\x=q_0(u\x_1+v\x_2)$.

Let 
\[U:=\left(q_1\frac{||\x_2||}{||\x_1||}\right)^{1/2}\andd
V:=\left(q_1\frac{||\x_1||}{||\x_2||}\right)^{1/2},\]
so that $UV=q_1$.  Then an easy application of the
pigeon-hole principle shows that one can find $(u,v)\in\Z^2-\{(0,0)\}$
with $q_1\mid L_1(u,v)$ and satisfying $|u|\le U$ and $|v|\le V$.  We
then deduce that 
\begin{eqnarray*}
||\x||&=&q_0||u\x_1+v\x_2||\\
&\le& q_0( U||\x_1||+V||\x_2||)\\
&=&2q_0\left(q_1||\x_1||.||\x_2||\right)^{1/2}\\
&\ll& q_0(q_1||\a_0||)^{1/2}\\
&=& q^{1/2}(q_0||\a_0||)^{1/2}\\
&\le &q^{1/2}(\alpha||\a_0||)^{1/2}\\
&=& q^{1/2}||\a||^{1/2}
\end{eqnarray*}
by (\ref{12a}). Since $\x$ must be non-zero we deduce that
\[m(Q;q)\ll (q||\a||)^{1/2}= q^{1/2}\widehat{m}(-Q^{\adj};q)^{1/2},\]
as required.

\section{Deduction of Theorem \ref{main}}\label{3}

In this section we will show how Theorem \ref{main} follows from
Theorem \ref{burgthm}.  Clearly it suffices to prove that
$\widehat{m}(Q;q)\ll_{\ep}q^{1/4+\ep}$ uniformly for any ternary form
$Q$ with $(\det(Q),q)=1$.

Our first task is to establish the following corollary to Theorem
\ref{burgthm}.
\begin{lemma}\label{exp}
For any $\delta>0$ there is a corresponding $\eta>0$ such that, if
$q>1$ and $C\subset \R^2$ is a convex
set contained in a disc $\{\x\in\R^2:\,||\x-\x_0||\le R\}$, then
\[\sum_{(x,y)\in C}\chi(Q(x,y))\ll_{\delta}R^{2-\eta}\]
for $R\ge q^{1/4+\delta}$, uniformly for every primitive
character $\chi$ modulo $q$ and for every binary
quadratic form $Q$ with $(\det(Q),q)=1$.
\end{lemma}

We prove this in three steps, beginning with the case in which
$q^{1/4+\delta}\le R\le q^{5/12}$. In this range we choose
\[r=3+[1/\delta],\;\;\; \eta=1/(r^2+4r)\andd \ep=\eta/4.\]
Then $r\ge 3$ and $1/4+\delta\ge 1/4+1/r$, so that
\[q\le R^{1/(1/4+\delta)}\le R^{4r/(r+4)}.\]
Thus Theorem \ref{burgthm} produces
\begin{eqnarray*}
\sum_{(x,y)\in C}\chi(Q(x,y))&\ll_{\ep,r}&R^{2-1/r}q^{(r+2)/(4r^2)+\ep}\\
&\ll_{\delta}&R^{2-1/r+(r+2)/r(r+4)+4\ep}\\
&=&R^{2-2/r(r+4)+\eta}\\
&=&R^{2-\eta}.
\end{eqnarray*}

Next, when $q^{5/12}\le R\ll q$, we cover $\R^2$ with disjoint squares
of side $q^{5/12}$ to obtain a partition of $C$ into $O(R^2q^{-5/6})$
convex subsets, each with diameter at most $q^{5/12}$.  On applying the
result above with $\delta=1/6$ we find that
\beql{star}
\sum_{(x,y)\in C}\chi(Q(x,y))\ll R^2q^{-5/6}
\big(q^{5/12}\big)^{2-\eta}\ll R^{2-5\eta/12}
\eeq
for some absolute constant $\eta>0$.

Finally we examine the case $R\gg q$. This time we cover $C$ with squares
of side $q$, and observe that 
\[\sum_{x,y\modd{q}}\chi(Q(x,y))=0.\]
(By multiplicativity it suffices to prove this when $q$ is
prime, in which case it is an easy exercise, relying on the fact that
$Q$ is nonsingular modulo $q$.) Since $C$ will be partitioned into
$O(R^2q^{-2})$ complete squares and $O(Rq^{-1})$ partial squares we may
use the result (\ref{star}) to conclude that
\[\sum_{(x,y)\in C}\chi(Q(x,y))\ll Rq^{-1}
q^{2-5\eta/12}\le R^{2-5\eta/12},\]
and the lemma follows.

We next estimate $\widehat{m}(Q;q)$ for binary forms $Q$.
\begin{lemma}\label{m2}
For any fixed $\delta>0$ we have
\[\widehat{m}(Q;q)\ll_{\delta}q^{1/4+\delta}\]
uniformly over odd square-free moduli $q$, and over binary forms $Q$
subject to $(\det(Q),q)=1$.
\end{lemma}

As already noted in the introduction, if we let $d$ run over
all divisors of $q$ then if $\sum_{d}\chi_d(m)>0$ we must
have $m\equiv t^2\modd{q}$ for some integer $t$.  It
follows that $\widehat{m}(Q;q)\le R$ provided that
\[\sum_{d\mid q}\sum_{||(x,y)||\le R}\chi_d(Q(x,y))>0.\]
The number of divisors of $q$ is  $O_{\ep}(q^{\ep})$, for any $\ep>0$.
Choosing $\ep=\eta/8$ it follows from Lemma \ref{exp} that if $R\ge
q^{1/4+\delta}$ then
\[\sum_{d\mid q,\,d>1}\sum_{||(x,y)||\le
  R}\chi_d(Q(x,y))\ll_{\delta}R^{2-\eta}q^{\ep}
\ll R^{2-\eta/2}.\]
On the other hand
\[\sum_{||(x,y)||\le R}1\gg R^2,\]
and Lemma \ref{m2} follows.

Finally we need to estimate $\widehat{B}_3(q)$ in terms of
$\widehat{B}_2(q)$.
\begin{lemma}\label{tra}
We have
\[\widehat{B}_3(q)\ll_{\ep}q^{\ep}\widehat{B}_2(q)\]
for any fixed $\ep>0$.
\end{lemma}
Once this is proved we may deduce from Lemma \ref{m2} that
$\widehat{B}_3(q)\ll_{\ep}q^{1/4+\ep}$ for any $\ep>0$, whence Lemma
\ref{key} yields $B_3^*(q)\ll_{\ep}q^{5/8+\ep}$. This is the result required for
Theorem \ref{main}.

To establish Lemma \ref{tra} we will find short vectors
\beql{vec}
(a_1,a_2),\,(a_3,a_4),\,(a_5,a_6)\in\Z^2
\eeq
such that the form
\beql{r1}
R(u,v):=Q(a_1u+a_2v,a_3u+a_4v,a_5u+a_6v)
\eeq
has $(\det(R),q)=1$.  We can then choose $u,v\ll \widehat{B}_2(q)$,
not both zero, such that $R(u,v)$ is a square modulo $q$, which will
produce a corresponding vector
\[\x=(a_1u+a_2v,a_3u+a_4v,a_5u+a_6v)\]
for which $Q(\x)$ is a square
modulo $q$, and with
\[||\x||\ll ||(u,v)||\max(a_1,\ldots,a_6).\]
If $\x$ were to vanish the three vectors (\ref{vec}) would all have
to be proportional.  But then the form (\ref{r1}) would have rank at
most 1,
so that $\det(R)=0$. This would contradict our assumption that
$(\det(R),q)=1$. It follows that we must have $\x\not=\mathbf{0}$.
Thus to complete the proof of Lemma \ref{tra} it will suffice to show
that we can choose the coefficients $a_1,\ldots,a_6$  to be of size
$O_{\ep}(q^{\ep})$. 

Define
\[\Delta(a_1,\ldots,a_6):=
\det\big(Q(a_1u+a_2v,a_3u+a_4v,a_5u+a_6v)\big).\]
This will be a sextic form in the 6 variables $a_1,\ldots,a_6$.  We
claim that for each prime factor $p$ of $q$ there is at least one
choice of $\a\in\Z^6$ such that $p\nmid\Delta(\a)$.  Since we can 
diagonalize $Q$ by a unimodular transformation over $\mathbb{F}_p$ a moment's
reflection shows that it is enough to verify the claim when $Q$ is a
diagonal form.  However the result is trivial in this case since
$p\nmid\det(Q)$. 

We can now call on the following lemma, which we will prove in a
moment.
\begin{lemma}\label{cop}
Let $\ep,\delta>0$ be given.  Suppose that
$F(x_1,\dots,x_n)\in\Z[x_1,\ldots,x_n]$ is a form of degree $d$,
and let $q\in\N$.  Assume that for every prime divisor $p$ of $q$
there is at least one $\a\in\Z^n$ such that $p\nmid F(\a)$.  Then 
\[\#\{\a\in\N^n:\,\max a_i\le A,\,(F(\a),q)=1\}\gg_{d,n,\ep}A^nq^{-\ep}\]
as soon as $A\ge q^{\delta}$ and $q\gg_{n,d}1$.
\end{lemma}
This result shows that we have at least one vector $\a$ of size
$||\a||\ll q^{\ep}$ such that $\Delta(\a)$ is coprime to $q$, which
suffices to complete the proof.

It remains to prove Lemma \ref{cop}. Define
\[N(e):=\#\{\a\modd{e}:\, e\mid F(\a)\}\]
for each $e\in\N$. Then $N(e)$ is multiplicative, and $N(p)<p^n$ for
$p\mid q$ by the hypothesis of the lemma. Moreover, when $N(p)<p^n$ the form
$F$ cannot vanish identically modulo $p$, whence $N(p)\ll_{d,n}
p^{n-1}$. It follows that 
\[N(e)\ll_{d,n,\eta}e^{n-1+\eta}\]
for $e\mid q$, for any fixed $\eta>0$.

We now consider 
\[N(e,A):=\#\{\a\in\N^n:\, \max a_i\le A,\, e\mid F(\a)\}.\]
The set $(0,A]^n$ contains $[A/e]^n$ disjoint cubes of side-length
$e$, and is included in a union of $(1+[A/e])^n$ such cubes. It
follows that 
\beql{both}
N(e,A)=\frac{A^n}{e^n}N(e)+O_{n}(A^{n-1}e^{1-n}N(e))
=\frac{A^n}{e^n}N(e)+O_{d,n,\eta}(A^{n-1}e^{\eta})
\eeq
when $e\mid q$ and $e\le A$. 
To handle larger values of $e$ we use a rather general result of
Browning and Heath-Brown \cite{bhb}.
For each $p_i\mid q$ let $V_i$ be the affine
variety over $\mathbb{F}_{p_i}$ given by $F=0$.
Since $F$ does not vanish identically modulo $p_i$ this has dimension
$n-1$.  We now apply \cite[Lemma 4]{bhb} with
$W=\mathbb{A}^n$, and $k_i=n-1$ for every index $i$.  Taking $e\mid q$ with
$e\ge A$
we find that there is a constant $C=C(d,n)$ such that
\[N(e,A)\ll C^{\omega(e)}(A^ne^{-1}+\omega(e)A^{n-1})\ll_{d,n,\eta}
e^{\eta}A^{n-1}\]
for any fixed $\eta>0$. It follows that if $e\ge A$ we will have
\begin{eqnarray*}
N(e,A)&=&\frac{A^n}{e^n}N(e)+O(A^nN(e)e^{-n})+O_{d,n,\eta}(A^{n-1}e^{\eta})\\
&=&\frac{A^n}{e^n}N(e)+O_{d,n,\eta}(A^ne^{-1+\eta})+O_{d,n,\eta}(A^{n-1}e^{\eta})\\
&=&\frac{A^n}{e^n}N(e)+O_{d,n,\eta}(A^{n-1}e^{\eta}),
\end{eqnarray*}
so that (\ref{both}) holds whether $e\le A$ or not.

We now examine
\begin{eqnarray*}
\lefteqn{\#\{\a\in\N^n:\,\max a_i\le A,\,(F(\a),q)=1\}}\\
&=&\sum_{e\mid q}\mu(e)N(e,A)\\
&=&\sum_{e\mid q}\mu(e)\frac{A^n}{e^n}N(e)+O_{d,n,\eta}
\left(\sum_{e\mid q} A^{n-1}e^{\eta}\right)\\
&=&A^n\prod_{p\mid q}(1-N(p)p^{-n})+O_{d,n,\eta}(A^{n-1}q^{2\eta}).
\end{eqnarray*}
Since $N(p)<p^n$ and $N(p)\le c_0 p^{n-1}$ for some constant $c_0$
depending only on $d$ and $n$, we may deduce that
\begin{eqnarray*}
\prod_{p\mid q}(1-N(p)p^{-n})&\ge& \prod_{\substack{p\mid q\\ p\le
    2c_0}}(1-(p^n-1)p^{-n})\prod_{\substack{p\mid
    q\\ p>2c_0}}(1-c_0p^{-1})\\
&\ge& \prod_{p\le 2c_0}p^{-n} \prod_{\substack{p\mid q\\ p>2c_0}}(1-p^{-1})^{2c_0}\\
&\gg_{d,n}&\left(\frac{\phi(q)}{q}\right)^{2c_0}\\
&\gg_{d,n,\eta}&q^{-\eta}.
\end{eqnarray*}
It follows that
\[\#\{\a\in\N^n:\,\max a_i\le A,\,(F(\a),q)=1\}
\ge c_2A^nq^{-\eta}-c_3A^{n-1}q^{2\eta}\]
for suitable positive constants $c_2$ and $c_3$ depending on $d$ and
$n$.  The lemma then follows on taking $\eta=\min(\ep,\delta/4)$.

\section{Proof of Theorem \ref{burgthm}}

For the proof we will write
\[\Sigma=\sum_{(x,y)\in C}\chi(Q(x,y))\]
for convenience. Let $N\in\N$ be a parameter to be chosen, satisfying
$N\le Rq^{-1/100}$, say, and set $S=[R/N]$.  We need to specify 
a ``good'' set of vectors
$\b{s}\in\N^2$, and this will require a further definition.  The form
$Q(X,Y)$ should be thought of as lying in $\left(\Z/q\Z\right)[X,Y]$,
  and we need an appropriate lift to $\Z[X,Y]$.  To achieve this we write
$Q(x_1,x_2)=Ax_1^2+Bx_1x_2+Cx_2^2$ and 
\beql{ld}
\Lambda=\{\b{v}\in\Z^3:\,\b{v}\equiv\lambda(A,B,C)\modd{q}\mbox{ for
  some }\lambda\in\Z\},
\eeq
and we let $(A^*,B^*,C^*)$ be a non-zero
vector in $\Lambda$ of minimal length. As there is a non-zero vector
$(A',B',C')\equiv(A,B,C)\modd{q}$ in $\Lambda$ with $|A'|,|B'|,|C'|\le
q/2$ we see that $q$ cannot divide $(A^*,B^*,C^*)$.  We now define
\[Q^*(X,Y)=A^*X^2+B^*XY+C^*Y^2.\]
Note that $\det(Q^*)\equiv\lambda^2\det(Q)\modd{q}$ for an appropriate
$\lambda$.  Since $q$ cannot
divide $\lambda$, and is
square-free and coprime to $\det(Q)$, we will have $q\nmid\det(Q^*)$.
In particular $Q^*$ is nonsingular, but there is no guarantee that
$(\det(Q^*),q)=1$. 
We can now take our set of ``good'' vectors $\b{s}$ to be
\[\clS=\{(s_1,s_2)\in\N^2:\, ||\bs||\le S,\, (Q(\bs),q)=1, 
Q^*(\b{s})\not=0\}.\]
There are $O(S)$ vectors for which $Q^*(\b{s})=0$, uniformly over all
non-zero forms $Q^*$.  Thus, according to Lemma \ref{cop} we have
\beql{lb}
\#\clS\gg_{\ep} S^2q^{-\ep},
\eeq
for $S\gg_{\ep} q^{\ep}$, for any fixed $\ep>0$.  

For any positive integer $n\le N$ we proceed to write
\[(\#\clS)\Sigma=\sum_{(s_1,s_2)\in\clS}\sum_{(x_1,x_2)\in Z^2}
\chi(Q(x_1+ns_1,x_2+ns_2))\mathds{1}_C(x_1+ns_1,x_2+ns_2),\]
where $\mathds{1}_C$ is the characteristic function for the set
$C$. It follows that
\[N(\#\clS)\Sigma=\sum_{(s_1,s_2)\in\clS}\sum_{(x_1,x_2)\in Z^2}
\sum_{n\in I}\chi(Q(x_1+ns_1,x_2+ns_2)),\]
where
\[I=\{n\le N:\,(x_1+ns_1,x_2+ns_2)\in C\}.\]
Since $C$ is convex, $I$ is an interval.  Moreover if $I$ is nonempty,
containing $\x+n\mathbf{s}$, then $||\x+n\mathbf{s}-\x_0||\le R$ and
$||n\mathbf{s}||\le NS\le R$, whence $||\x-\x_0||\le 2R$. We therefore
deduce, via (\ref{lb}), that
\[\Sigma\ll_{\ep}N^{-1}S^{-2}q^{\ep}\sum_{\mathbf{s}\in\clS}
\sum_{\substack{\x\in Z^2\\ ||\x-\x_0||\le 2R}}
\max_{I\subseteq(0,N]}\left|\sum_{n\in I}\chi(Q(x_1+ns_1,x_2+ns_2))\right|.\]

If the reader compares this with the corresponding stage in the
argument of Chang \cite{ch}, see \cite[(4.3)]{ch} for example, then it will be
observed that Chang has a product $st$ in place of our variable
$n$. Indeed our method is slightly different from Chang's, 
requires one variable fewer, and does not use an argument
corresponding to \cite[Lemma 3]{ch}.

To proceed further we use the readily verified identity
\[Q(x_1+ns_1,x_2+ns_2)=Q(\bs)\widetilde{Q}\big(n+a(\bs,\x),b(\bs,\x)\big),\]
where $\widetilde{Q}(x_1,x_2):=x_1^2+Bx_1x_2+ACx_2^2$, and
\beql{abd}
a(\bs,\x)=\frac{Ax_1s_1+Bx_1s_2+Cx_2s_2}{Q(s_1,s_2)},\;\;\;
b(\bs,\x)=\frac{x_2s_1-x_1s_2}{Q(s_1,s_2)}.
\eeq
Here the fractions are to be interpreted in the ring $\Z/q\Z$, the
denominators $Q(s_1,s_2)$ being units by our choice of the set $\clS$.
We now write
\[N(a,b)=\#\left\{(\mathbf{s},\x)\in\clS\times\Z^2:\,||\x-\x_0||\le 2R,\,
a(\bs,\x)=a,\,b(\bs,\x)=b\right\},\]
whence
\[\Sigma\ll_{\ep}N^{-1}S^{-2}q^{\ep}\sum_{a,b\modd{q}}N(a,b)
\max_{I\subseteq(0,N]}\left|\sum_{n\in I}\chi(\widetilde{Q}(n+a,b))\right|.\]

We must now consider the mean square of $N(a,b)$, for which we will
prove the following bound.
\begin{lemma}\label{last}
For any fixed $\ep>0$ we have
\[\sum_{a,b\modd{q}}N(a,b)^2\ll_{\ep}q^{\ep}R^2S^2(1+RSq^{-1/2}+R^2S^2q^{-4/3}).\]
\end{lemma}
This will be established in the next section.

We also have the trivial bound
\[\sum_{a,b\modd{q}}N(a,b)\le\#\{(\mathbf{s},\x)\in\clS\times\Z^2:\,
||\x-\x_0||\le 2R\} \ll R^2S^2,\]
whence H\"{o}lder's inequality yields
\begin{eqnarray}\label{hol}
\Sigma^{2r}
&\ll_{\ep,r}&\left(N^{-1}S^{-2}q^{\ep}\right)^{2r}
\left\{\sum_{a,b\modd{q}}N(a,b)\right\}^{2r-2}
\left\{\sum_{a,b\modd{q}}N(a,b)^2\right\}\nonumber\\
&&\hspace{2cm}\mbox{}\times
\sum_{a,b\modd{q}}\max_{I\subseteq(0,N]}
\left|\sum_{n\in I}\chi(\widetilde{Q}(n+a,b))\right|^{2r}\nonumber\\
&\ll_{\ep,r}&
N^{-2r}R^{4r-2}S^{-2}q^{\ep}(1+RSq^{-1/2}+R^2S^2q^{-4/3})\nonumber\\
&&\hspace{2cm}\mbox{}\times\sum_{a,b\modd{q}}
\max_{I\subseteq(0,N]}
\left|\sum_{n\in I}\chi(\widetilde{Q}(n+a,b))\right|^{2r}\nonumber\\
&\ll_{\ep,r}&
N^{2-2r}R^{4r-4}q^{\ep}(1+R^2N^{-1}q^{-1/2}+R^4N^{-2}q^{-4/3})\nonumber\\
&&\hspace{2cm}\mbox{}\times
\sum_{a,b\modd{q}}\max_{I\subseteq(0,N]}
\left|\sum_{n\in I}\chi(\widetilde{Q}(n+a,b))\right|^{2r},
\end{eqnarray}
on employing our convention concerning the values taken by $\ep$.

We are therefore led to consider sums of the form
\[S(q;H):=\sum_{a,b\modd{q}}
\left|\sum_{n\le H}\chi(\widetilde{Q}(n+a,b))\right|^{2r}.\]
To estimate these we expand to obtain
\[S(q;H)=\sum_{n_1,\ldots,n_{2r}\le H}\Sigma(q;\b{n})\]
with
\[\Sigma(q;\b{n})=\sum_{a,b\modd{q}}\chi(F_+(a,b;\b{n}))
\overline{\chi}(F_-(a,b;\b{n}))\]
and
\[F_+(X,Y;\b{n})=\prod_{i=1}^r\widetilde{Q}(n_i+X,Y),\;\;\;
F_-(X,Y;\b{n})=\prod_{i=r+1}^{2r}\widetilde{Q}(n_i+X,Y).\]
The sums $\Sigma(q;\b{n})$ have a standard multiplicative
property.  If $q=uv$, say, then $u$ and $v$ will be coprime and
square-free, and we can write $\chi=\chi_u\chi_v$ for suitable
primitive characters to moduli $u$ and $v$ respectively.  We will then
have
\beql{mult}
\Sigma(q;\b{n})=\Sigma(u;\b{n})\Sigma(v;\b{n}).
\eeq
It therefore suffices to understand $\Sigma(q;\b{n})$ when $q$ is
prime, for which we have the following result.
\begin{lemma}\label{wl}
Let $p$ be an odd prime not dividing $\det(\widetilde{Q})$, and let
$\chi$ be a non-principal character to
modulus $p$.  Write
\[\Delta_i=\prod_{\substack{1\le j\le 2r\\ j\not=i}}(n_j-n_i)\]
and
\[\Delta=\hcf(\Delta_1,\ldots,\Delta_{2r}).\]
Then
\[|\Sigma(p;\b{n})|\le 4r^2p(p,\Delta).\]
\end{lemma}

We will prove this in Section \ref{ltt}.  By summing over the
$(2r)$-tuples $\b{n}$ we are then able to establish the following
bound for $S(q;H)$.
\begin{lemma}\label{SqH}
For any $\ep>0$ and $r\in\N$ we have
\[S(q,H)\ll_{\ep,r}(qH)^{\ep}(qH^{2r}+q^2H^r).\]
\end{lemma}
This will be proved in Section \ref{ls}

Having established this there is a standard procedure to insert a
maximum over subintervals of $(0,N]$, which goes back to
Rademacher \cite{Rad} and Menchov \cite{Men}. We do not repeat the
details, but instead refer the reader to Gallager and Montgomery
\cite[Section3]{GM} or Heath-Brown \cite[Section 2]{bk}.  The outcome
is the following result.

\begin{lemma}\label{SqHm}
For any $\ep>0$ and $r\in\N$ we have
\[\sum_{a,b\modd{q}}\max_{I\subseteq(0,N]}
\left|\sum_{n\in I}\chi(\widetilde{Q}(n+a,b))\right|^{2r}
\ll_{\ep,r}(qN)^{\ep}(qN^{2r}+q^2N^r).\]
\end{lemma}

We are now ready to complete the proof of Theorem \ref{burgthm}.  We
insert the bound of Lemma \ref{SqHm} into (\ref{hol}), to give
\[\Sigma^{2r}\ll_{\ep,r}N^{2-2r}R^{4r-4}q^{\ep}
(1+R^2N^{-1}q^{-1/2}+R^4N^{-2}q^{-4/3}).(qN)^{\ep}(qN^{2r}+q^2N^r).\]
In order to balance the final two terms we choose $N=[q^{1/r}]$, which
satisfies our constraint $N\le Rq^{-1/100}$ provided that 
$R\ge q^{1/4+1/2r}$ and
$r\ge 3$. On re-defining $\ep$ we then
find that
\begin{eqnarray*}
\Sigma^{2r}&\ll_{\ep,r}&q^{\ep}N^{2-2r}R^{4r-4}
(1+R^2N^{-1}q^{-1/2}+R^4N^{-2}q^{-4/3}).qN^{2r}\\
&\ll_{\ep,r}&q^{1/2+1/r+\ep}R^{4r-2}(R^{-2}q^{1/2+1/r}+1+R^2q^{-5/6-1/r}),
\end{eqnarray*}
and the theorem follows.

\section{Proof of Lemma \ref{last}}
We now prove Lemma \ref{last}.  
In view of the definitions
(\ref{abd}) we have the identity
\[(Ab(\bs,\x)X-a(\bs,\x)Y)(s_2X-s_1Y)=s_2b(\bs,\x)Q(X,Y)-(x_2X-x_1Y)Y\]
in $\big(\Z/q\Z\big)[X,Y]$.
Thus if $a(\bs,\x)=a(\bs',\x')=a$ and $b(\bs,\x)=b(\bs',\x')=b$ then
\begin{eqnarray*}
\lefteqn{(AbX-aY)(s_2X-s_1Y)(s_2'X-s_1'Y)}\hspace{2cm}\\
&=&\big(s_2bQ(X,Y)-(x_2X-x_1Y)Y\big)(s_2'X-s_1'Y),
\end{eqnarray*}
and also
\begin{eqnarray*}
\lefteqn{(AbX-aY)(s_2'X-s_1'Y)(s_2X-s_1Y)}\hspace{2cm}\\
&=&\big(s_2'bQ(X,Y)-(x_2'X-x_1'Y)Y\big)(s_2X-s_1Y).
\end{eqnarray*}
Thus by subtraction we deduce that
\begin{eqnarray*}
\lefteqn{Y\{(x_2X-x_1Y)(s_2'X-s_1'Y)-(x_2'X-x_1'Y)(s_2X-s_1Y)\}}\hspace{4cm}\\
&=& b(s_2's_1-s_1's_2)YQ(X,Y),
\end{eqnarray*}
still in $\big(\Z/q\Z\big)[X,Y]$.

We then deduce that
\beql{lam}
\big(x_2s_2'-x_2's_2,\;x_2's_1+x_1's_2-x_2s_1'-x_1s_2',\;
x_1s_1'-x_1's_1\big)\in\Lambda,
\eeq
with $\Lambda$ given by (\ref{ld}).  It follows that 
\[\sum_{a,b}N(a,b)^2\le\sum_{\bs,\bs'\in\clS}N_1(\bs,\bs'),\]
where $N_1(\bs,\bs')$ counts pairs of vectors $(\x,\x')$ each lying in the disc
$||\x-\x_0||\le 2R$, such that (\ref{lam}) holds.  Now suppose that
$(\x_1,\x_1')$ is a pair counted by $N_1(\bs,\bs')$.  For any other
such pair we write $\x=\x_1+\b{u}$ and $\x'=\x_1'+\b{u}'$ whence, by
subtraction, we find firstly that $||\b{u}||,||\b{u}'||\le 4R$, and
secondly that
\beql{lam1}
(u_2s_2'-u_2's_2,\,u_2's_1+u_1's_2-u_2s_1'-u_1s_2',\,
u_1s_1'-u_1's_1)\in\Lambda.
\eeq
Thus $N_1(\bs,\bs')\le N_2(\bs,\bs')$, where $N_2(\bs,\bs')$ counts
pairs of vectors $\b{u},\b{u}'$ satisfying (\ref{lam1}), and having
length at most $4R$. 

We have already chosen $(A^*,B^*,C^*)=\b{v}_1$, say, as the shortest
vector in $\Lambda$.  As in the proof of Davenport \cite[Lemma5]{16},
we can then construct a basis $\b{v}_1,\b{v}_2,\b{v}_3$ for $\Lambda$,
such that if
$\b{v}=\lambda_1\b{v}_1+\lambda_2\b{v}_2+\lambda_3\b{v}_3$, then
$\lambda_i\ll||\b{v}||/||\b{v}_i||$ for $i=1,2,3$. Moreover we will
have
\[||\b{v}_1||\le||\b{v}_2||\le||\b{v}_3||\]
and
\[||\b{v}_1||.||\b{v}_2||.||\b{v}_3||\ge\det(\Lambda).\]
In our case we have $\det(\Lambda)=q^2$, whence
\[||\b{v}_2||.||\b{v}_3||\ge q^{4/3}.\]
Moreover one sees from the definition of $\Lambda$ that $q\mid
\b{v}_1\wedge\b{v}_2$, and since the vectors $\b{v}_1$ and
$\b{v}_2$ are not proportional it follows that
\[q\le||\b{v}_1\wedge\b{v}_2||\le||\b{v}_1||.||\b{v}_2||\le
||\b{v}_2||^2.\]

The vector 
\[\b{v}=(u_2s_2'-u_2's_2,\,u_2's_1+u_1's_2-u_2s_1'-u_1s_2',\,
u_1s_1'-u_1's_1)\]
has length at most $32RS$ so that the corresponding coefficients
satisfy 
\[\lambda_2\ll \frac{RS}{||\b{v}_2||}\;\;\;\mbox{and}\;\;\;
\lambda_3\ll \frac{RS}{||\b{v}_3||}.\]
If we break the available vectors counted by $N_2(\bs,\bs')$ into
subsets according to the values of $\lambda_2$ and $\lambda_3$, then
the number of such subsets will be
\[\ll \left(1+\frac{RS}{||\b{v}_2||}\right)
\left(1+\frac{RS}{||\b{v}_3||}\right)\ll
1+\frac{RS}{||\b{v}_2||}+\frac{R^2S^2}{||\b{v}_2||.||\b{v}_3||}
\ll 1+\frac{RS}{q^{1/2}}+\frac{R^2S^2}{q^{4/3}}.\]
If $(\b{u}_1,\b{u}_1')$ and $(\b{u}_2,\b{u}_2')$ are two pairs
belonging to the same subset, and we write $\b{u}=\b{u}_1-\b{u}_2$ and
$\b{u}'=\b{u}_1'-\b{u}_2'$, then  
\beql{vf}
\b{v}=(u_2s_2'-u_2's_2,\,u_2's_1+u_1's_2-u_2s_1'-u_1s_2',\,
u_1s_1'-u_1's_1)
\eeq
will be a multiple of $\b{v}_1=(A^*,B^*,C^*)$, and we will have
$||\b{u}||,||\b{u}'||\le 8R$.

We therefore conclude that
\[N_2(\bs,\bs')\ll(1+RSq^{-1/2}+R^2S^2q^{-4/3})N_3(\bs,\bs')\]
where $N_3(\bs,\bs')$ counts pairs of vectors $\b{u},\b{u}'$ having
length at most $8R$, and for which the vector (\ref{vf}) is an integer 
multiple of $(A^*,B^*,C^*)$.  The quadratic form
corresponding to $\b{v}$ is
\begin{eqnarray*}
\lefteqn{(u_2s_2'-u_2's_2)X^2+(u_2's_1+u_1's_2-u_2s_1'-u_1s_2')XY
+(u_1s_1'-u_1's_1)Y^2}\\
&=&(u_2X-u_1Y)(s_2'X-s_1'Y)-(u_2'X-u_1'Y)(s_2X-s_1Y).
\end{eqnarray*}
We therefore conclude that
\beql{div}
Q^*(X,Y)\mid(u_2X-u_1Y)(s_2'X-s_1'Y)-(u_2'X-u_1'Y)(s_2X-s_1Y).
\eeq
Thus, to complete the proof of Lemma \ref{last} it 
suffices to show that 
\[\#\{(\b{u},\b{u}',\b{s},\b{s}')\in\Z^2\times\Z^2\times\clS\times\clS:\,
||\b{u}||,||\b{u}'||\le 8R,\mbox{ (\ref{div}) holds}\}\]
\beql{lt}
\ll_{\ep}q^{\ep}R^2S^2.
\eeq

Given two binary quadratic forms $Q_1$ and $Q_2$ one may define a
covariant $C(Q_1,Q_2)$ as the discriminant of the binary form
$D(\alpha,\beta)=\det(\alpha Q_1+\beta Q_2)$. One readily confirms
that $C(Q_1,Q_2)=C(Q_1+\lambda Q_2,Q_2)$ for any constant $\lambda$,
and moreover that
\[C\big((u_2X-u_1Y)(v_2X-v_1Y),Q\big)=Q(u_1,u_2)Q(v_1,v_2).\]
Taking $Q_1=(u_2X-u_1Y)(s_2'X-s_1'Y)$ and $Q_2=Q^*$ we deduce that
\beql{=}
Q^*(u_1,u_2)Q^*(s_1',s_2')=Q^*(u_1',u_2')Q^*(s_1,s_2).
\eeq

In defining the set $\clS$ we arranged that $Q^*(\b{s})\not=0$. If
$Q^*(\b{u})\not=0$ then $Q^*(u_1,u_2)Q^*(s_1',s_2')$ has
$O_{\ep}(q^{\ep})$ divisors, since
\[|Q^*(u_1,u_2)Q^*(s_1',s_2')|\ll
\max(|A^*|,|B^*|,|C^*|)^2||\b{u}||^2||\b{s}||^2\ll q^2 R^2S^2\ll
q^6.\]
Moreover, when $d\not=0$ the equation $Q^*(u_1',u_2')=d$ will have
$\ll_{\ep}(qR)^{\ep}\ll_{\ep} q^{\ep}$ solutions $\b{u}'$ with
$||\b{u}'||\le 8R$, by Theorem 13 of Heath-Brown \cite{aa} for
example.  (Here we use crucially the fact that $Q^*$ is nonsingular.)
Similarly $Q^*(s_1,s_2)=d'$ will have $O_{\ep}(q^{\ep})$ solutions for
any $d'\not=0$.  It then follows that the contribution 
arising from 4-tuples $(\b{u},\b{u}',\b{s},\b{s}')$ in which
$Q^*(\b{u})\not=0$ will be $O_{\ep}(q^{\ep}R^2S^2)$, which is
satisfactory for (\ref{lt}).

It remains to deal with the case in which $Q^*(\b{u})=0$.  In view of
(\ref{=}) we will then have $Q^*(\b{u}')=0$, since $Q^*(\b{s})$ and
$Q^*(\b{s}')$ are non-zero.  We now claim that either
$\b{u}=\b{u}'=\b{0}$, or $Q^*(X,Y)$ factors over $\Z$
into linear factors $L_1(X,Y)$ and $L_2(X,Y)$ such that $L_1(X,Y)$ divides both
$u_2X-u_1Y$ and $u_2'X-u_1'Y$. To see this, suppose that
$\b{u}\not=\b{0}$, say. Then we must have $L_1(X,Y)\mid u_2X-u_1Y$ for
some integral linear factor of $Q^*$. It would then follow from
(\ref{div}) that $L_1(X,Y)\mid u_2'X-u_1'Y$, since $Q^*(\b{s})\not=0$.  The
claim then follows.

Clearly the contribution to (\ref{lt}) arising from the case
$\b{u}=\b{u}'=\b{0}$ is $O(S^4)=O(R^2S^2)$, which is
satisfactory, so it remains to consider the case
in which 
\[u_2X-u_1Y=kL_1(X,Y)\;\;\; \mbox{and}\;\;\; u_2'X-u_1'Y=k'L_1(X,Y)\]
with integers $k,k'$ such that $|k|,|k'|\le 8R$.  We then have
\[k(s_2'X-s_1'Y)\equiv k'(s_2X-s_1Y)\modd{L_2(X,Y)},\]
by (\ref{div}).
If $L_2(X,Y)=aX-bY$, say then we must have
\[k(s_2'b-s_1'a)=k'(s_2b-s_1a).\]
Moreover $s_2'b-s_1'a$ and $s_2b-s_1a$ are non-zero, since
$Q^*(\b{s})$ and $Q^*(\b{s}')$ do not vanish. If $k'=0$ then $k=0$,
which would put us in the case $\b{u}=\b{u}'=\b{0}$ which has already
been dealt with.  Since at most one of $a$ or $b$ can vanish we may
suppose that $b$ say, is non-zero. There are then $O(RS^3)$
possibilities for $s_1',s_1,s_2$ and $k'$, and the number of divisors
$k$ of $k'(s_2b-s_1a)$ will be $O_{\ep}((qS)^{\ep})$, since
$|a|,|b|\ll\max(|A^*|,|B^*|,|C^*|)\ll q$. The complementary divisor to $k$
is then $s_2'b-s_1'a$, which determines $s_2'$.  We therefore conclude
that the corresponding contribution to (\ref{lt}) is
$O_{\ep}(q^{\ep}R^2S^2)$, since $S\le R\le q$.  This completes the
proof of (\ref{lt}), and hence of Lemma \ref{last}.

\section{Proof of Lemma \ref{wl}}\label{ltt}
Our proof of Lemma \ref{wl} is inspired by the viewpoint taken by 
Chang \cite{ch}. We first consider the case in which
$\widetilde{Q}(X,Y)=X^2+BXY+ACY^2$ factors modulo $p$.  In this case
we may replace $\widetilde{Q}(X,Y)$ by $(X+\lambda Y)(X+\mu Y)$ say,
where $p\nmid\lambda-\mu$ since $p\nmid\det(\widetilde{Q})$. Then 
$\widetilde{Q}(n+a,b)=(n+a')(n+b')$ where $a'=a+\lambda b$ and
$b'=a+\mu b$ are independent of $n$.  Moreover $(a',b')$ runs over
$\mathbb{F}_p^2$ as $(a,b)$ does.  It follows that
\[\Sigma(p;\b{n})=\sum_{a,b\modd{p}}\chi(G_+(a,b;\b{n}))
\overline{\chi}(G_-(a,b;\b{n}))\]
with
\[G_+(X,Y;\b{n})=\prod_{i=1}^r(n_i+X)(n_i+Y),\;\;\;
G_-(X,Y;\b{n})=\prod_{i=r+1}^{2r}(n_i+X)(n_i+Y).\]
We then see that
\[\Sigma(p;\b{n})=\Sigma_1(p;\b{n})^2\]
with
\[\Sigma_1(p;\b{n})=\sum_{a\modd{p}}\chi(H_+(a;\b{n}))
\overline{\chi}(H_-(a;\b{n}))\]
and
\[H_+(X;\b{n})=\prod_{i=1}^r(n_i+X),\;\;\;
H_-(X;\b{n})=\prod_{i=r+1}^{2r}(n_i+X).\]

The sum $\Sigma_1(p;\b{n})$ occurs in the work of Burgess \cite[Lemma
1]{bprim}, from which one readily sees that
\beql{bw}
|\Sigma_1(p;\b{n})|\le 2r\sqrt{p}
\eeq
unless every linear factor of the polynomial
$H_+(X;\b{n})H_-(X;\b{n})$ has multiplicity two or more, modulo $p$.
In the exceptional case we have $p\mid\Delta_i$ for every $i$, 
whence $p\mid\Delta$.  We
deduce that (\ref{bw}) holds whenever $p\nmid\Delta$.  In the
remaining case we have a trivial bound $|\Sigma_1(p;\b{n})|\le p$, so
that
\[|\Sigma_1(p;\b{n})|\le 2rp^{1/2}(p,\Delta)^{1/2}\]
whether or not $p\nmid\Delta$.  We therefore conclude that
\[|\Sigma(p;\b{n})|\le 4r^2p(p,\Delta)\]
whenever $\widetilde{Q}$ factors modulo $p$.  This is satisfactory for
Lemma \ref{wl}.

We turn now to the case in which $\widetilde {Q}$ is irreducible over
$\mathbb{F}_p$. It will be typographically convenient to write
$F$ for the field $\mathbb{F}_{p^2}$. In the case under
consideration, there is a factorization
$\widetilde{Q}(X,Y)=(X+\lambda Y)(X+\lambda' Y)$ say, over $F$
with $\lambda$ and $\lambda'$ being conjugates in
$F/\mathbb{F}_p$. We may now define a function $\psi$
from $F$ to $\C$ by setting
\[\psi(a+\lambda b)=\chi\big((a+\lambda b)(a+\lambda' b)\big)=
\chi(\widetilde{Q}(a,b)).\]
One easily sees that this is a non-trivial multiplicative character on 
$F$, and that
\[\Sigma(p;\b{n})=\sum_{\alpha\in F}\psi(H_+(\alpha;\b{n}))
\overline{\psi}(H_-(\alpha;\b{n})).\]
Burgess' proof of (\ref{bw}), based on Weil's ``Riemann Hypothesis''
for curves over arbitrary finite fields, immediately extends to
$\Sigma(p;\b{n})$, and shows that
\[|\Sigma(p;\b{n})|\le 2r\sqrt{\# F}=2rp\]
unless every linear factor of the polynomial
$H_+(X;\b{n})H_-(X;\b{n})$ has multiplicity two or more, modulo $p$.
In the alternative case we have $p\mid\Delta$, in the notation of the lemma,
and we deduce that
\[|\Sigma(p;\b{n})|\le p(p,\Delta),\]
in view of the trivial bound $|\Sigma(p;\b{n})|\le p^2$. As above, these
bounds are satisfactory for
Lemma \ref{wl}.

\section{Proof of Lemma \ref{SqH}}\label{ls}

It follows from Lemma \ref{wl}, along with the multiplicative relation
(\ref{mult}) that
\[\Sigma(q;\b{n})\le (4r^2)^{\omega(q)}q(q,\Delta)\ll_{\ep,r}
q^{1+\ep}(q,\Delta).\]  Thus to prove Lemma \ref{SqH} it will be
enough to show that
\[\sum_{n_1,\ldots,n_{2r}\le H}(q,\Delta)\ll_{\ep,r}
(qH)^{\ep}(H^{2r}+qH^r).\]
Indeed, we have
\[\sum_{n_1,\ldots,n_{2r}\le H}(q,\Delta)\le
\sum_{k\mid q}k\#\{\b{n}\in\N^{2r}\cap(0,H]^{2r}:\, k\mid\Delta\},\]
so that it suffices to establish the estimate
\beql{bb}
\#\{\b{n}\in\N^{2r}\cap(0,H]^{2r}:\, k\mid\Delta\}
\ll_{\ep,r} (kH)^{\ep}(H^{2r}k^{-1}+H^r).
\eeq

We first consider vectors $\b{n}$ for
$\Delta_1=\ldots=\Delta_{2r}=0$. Then if $\nu\in\N$ and there is any
index $i$ such that $n_i=\nu$, there must be at least two such
indices.  It follows that the set $\{n_1,\ldots,n_{2r}\}$ contains at
most $r$ distinct elements. There are at most $H^r$ choices for these
elements, $\nu_1,\ldots,\nu_s$ say, with $1\le s\le r$. Once the
$\nu_j$ have been chosen there are (at most) $s$ choices for each
$n_i$.  It follows that there are $O_r(H^r)$ vectors $\b{n}$ for
which $\Delta_1=\ldots=\Delta_{2r}=0$.  This is satisfactory for
(\ref{bb}).

In the remaining case we have $\Delta_j\not=0$ for some index $j$, and
\begin{eqnarray*}
\lefteqn{\#\{\b{n}\in\N^{2r}\cap(0,H]^{2r}:\,
  k\mid\Delta,\,\Delta\not=0\}}
\hspace{1cm}\\
&\le& \sum_{j=1}^{2r}
\#\{\b{n}\in\N^{2r}\cap(0,H]^{2r}:\, k\mid\Delta_j,\,\Delta_j\not=0\}.
\end{eqnarray*}
However $|\Delta_j|\le H^{2r-1}$, so that there are at most
$2H^{2r-1}k^{-1}$ possibilities for $\Delta_j$. For each such choice of
$\Delta_j$ there are at most $2d(|\Delta_j|)\ll_{\ep,r}H^{\ep}$
possibilities for each of its divisors $n_i-n_j$. Thus, taking account
of the $O(H)$ possibilities for $n_j$ itself, we find that
\[\#\{\b{n}\in\N^{2r}\cap(0,H]^{2r}:\, k\mid\Delta_j,\,\Delta_j\not=0\}
\ll_{\ep,r} H^{2r-1}k^{-1}(H^{\ep})^{2r-1}H.\]
After replacing $\ep$ by $\ep/(2r-1)$ we see that this is 
$O_{\ep,r}(H^{2r+\ep}k^{-1})$.  Since this is satisfactory for
(\ref{bb}) the proof of Lemma \ref{SqH} is complete.

\section{Acknowledgement}
This research was supported by EPSRC grant EP/K021132X/1.

\bigskip
\bigskip

Mathematical Institute,

Radcliffe Observatory Quarter,

Woodstock Road,

Oxford

OX2 6GG

UK
\bigskip

{\tt rhb@maths.ox.ac.uk}

\end{document}